\documentclass[11pt]{article}

\usepackage{amscd,amsmath, amssymb, fancyhdr, graphicx, url}

\usepackage[backref=page]{hyperref}

\renewcommand*{\backrefalt}[4]{%
	\ifcase #1 (Not cited.)%
	\or        (Cited on page~#2.)%
	\else      (Cited on pages~#2.)%
	\fi}

\hypersetup{
	colorlinks   = true,
	citecolor    = magenta,
	linkcolor    =blue,
	urlcolor     =magenta	
}

\numberwithin{equation}{section}

\newcommand{\version}{version 1.0.2,\ \ Nov 7, 2023}

\makeatletter
\def\x@arrow{\DOTSB\Relbar}
\def\xlongrightarrowfill@{\arrowfill@\relbar\relbar\longrightarrow}
\newcommand{\xlongrightarrow}[2][]{%
        \ext@arrow 0099\xlongrightarrowfill@{#1}{#2}}
\makeatother

\def\eqref#1{(\ref{#1})}

\newcommand{\arrow}{{\:\longrightarrow\:}}
\newcommand{\Z}{{\Bbb Z}}
\def\C{{\Bbb C}}
\def\P{{\Bbb P}}

\newcommand{\R}{{\Bbb R}}

\def\1{\sqrt{-1}\:}
\newcommand{\restrict}[1]{{\left|_{{\phantom{|}\!\!}_{#1}}\right.}}
\newcommand{\cntrct}                
{\hspace{2pt}\raisebox{1pt}{\text{$\lrcorner$}}\hspace{2pt}}

\newcommand{\calo}{{\cal O}}

\renewcommand{\tilde}{\widetilde}
\renewcommand{\bar}{\overline}
\renewcommand{\phi}{\varphi}
\renewcommand{\epsilon}{\varepsilon}
\renewcommand{\geq}{\geqslant}
\renewcommand{\leq}{\leqslant}

\newcommand{\Tot}{\operatorname{Tot}}

\newcommand{\Sym}{\operatorname{Sym}}
\newcommand{\Pic}{\operatorname{Pic}}

\newcommand{\Supp}{\operatorname{Supp}}

\newcommand{\Diff}{\operatorname{\sf Diff}}

\newcommand{\proof}{\noindent{\bf Proof:\ }}
\newcommand{\pstep}{{\bf Proof. Step 1:\ }}

\newcounter{Mycounter}[section]
\newcounter{lemma}[section]
\renewcommand{\thelemma}{{Lemma \thesection.\arabic{lemma}}}
\newcommand{\lemma}{%
    \setcounter{lemma}{\value{Mycounter}}
    \refstepcounter{lemma}
    \stepcounter{Mycounter}
    {\noindent \bf \thelemma:\ }}

\newcounter{claim}[section]
\renewcommand{\theclaim}{{Claim \thesection.\arabic{claim}}}
\newcommand{\claim}{%
    \setcounter{claim}{\value{Mycounter}}
    \refstepcounter{claim}
    \stepcounter{Mycounter}
    {\noindent \bf \theclaim:\ }}

\newcounter{sublemma}[section]
\newcounter{corollary}[section]
\renewcommand{\thecorollary}{{Corollary \thesection.\arabic{corollary}}}
\newcommand{\corollary}{%
    \setcounter{corollary}{\value{Mycounter}}
    \refstepcounter{corollary}
    \stepcounter{Mycounter}
    {\noindent \bf \thecorollary:\ }}

\newcounter{theorem}[section]
\renewcommand{\thetheorem}{{Theorem \thesection.\arabic{theorem}}}
\newcommand{\theorem}{%
    \setcounter{theorem}{\value{Mycounter}}
    \refstepcounter{theorem}
    \stepcounter{Mycounter}
    {\noindent \bf \thetheorem:\ }}

\newcounter{conjecture}[section]
\newcounter{proposition}[section]
\renewcommand{\theproposition}
      {{Proposition \thesection.\arabic{proposition}}}
\newcommand{\proposition}{%
    \setcounter{proposition}{\value{Mycounter}}
    \refstepcounter{proposition}
    \stepcounter{Mycounter}
    {\noindent \bf \theproposition:\ }}

\newcounter{definition}[section]
\renewcommand{\thedefinition}
      {{Definition~\thesection.\arabic{definition}}}
\newcommand{\definition}{%
    \setcounter{definition}{\value{Mycounter}}
    \refstepcounter{definition}
    \stepcounter{Mycounter}
    {\noindent \bf \thedefinition:\ }}

\newcounter{example}[section]
\newcounter{remark}[section]
\renewcommand{\theremark}{{Remark \thesection.\arabic{remark}}}
\newcommand{\remark}{%
    \setcounter{remark}{\value{Mycounter}}
    \refstepcounter{remark}
    \stepcounter{Mycounter}
    {\noindent \bf \theremark:\ }}

\newcounter{problem}[section]
\newcounter{question}[section]
\makeatletter

\@addtoreset{equation}{section} \@addtoreset{footnote}{section}
\makeatother

\def\blacksquare{\hbox{\vrule width 5pt height 5pt depth 0pt}}
\def\endproof{\;\blacksquare}

\begin{document}

\begin{center}
{\LARGE\bf
Normal form  of bimeromorphically contractible
holomorphic Lagrangian submanifolds\\[3mm]
}

Ekaterina Amerik\footnote{Both authors acknowledge support of 
HSE University basic research program; also partially supported by ANR (France) project FANOHK} 
Misha Verbitsky\footnote{Partially supported by 
FAPERJ E-26/202.912/2018 and CNPq - Process 310952/2021-2. \\

{\bf Keywords:} hyperk\"ahler manifold, Lagrangian submanifold, 
normal form, holomorphically symplectic manifold

{\bf 2010 Mathematics Subject
Classification: 53C26} }

\end{center}

{\small \hspace{0.15\linewidth}
\begin{minipage}[t]{0.7\linewidth}
{\bf Abstract} \\
Let $M$ be a holomorphically symplectic complex manifold,
not necessarily compact or quasiprojective,
and $X \subset M$ a compact Lagrangian submanifold. 
We construct a deformation to the normal cone,
showing that a neighbourhood of $X$ can be deformed
to its neighbourhood in $T^* X$. This is used to 
study Lagrangian submanifolds which can be bimeromorphically contracted
to a point. We prove that such submanifolds are
biholomorphic to $\C P^n$, and show that a 
certain neighbourhood of $X$ is symplectically
biholomorphic to a neighbourhood of the zero section
of its cotangent bundle. This gives a holomorphic 
version of the Weinstein's normal neighbourhood theorem.
\end{minipage}
}

\tableofcontents


\section{Introduction}

\subsection{Holomorphically symplectic Moser's lemma 
and the deformation to the normal cone}

A smooth complex manifold $M$ is {\bf holomorphically
  symplectic} if it carries a closed holomorphic $2$-form
$\Omega$, non-degenerate at each point. Clearly, in this
case $M$ is of even dimension $2n$ and the canonical
bundle of $M$ is trivial. A subvariety $X\subset M$ is
called {\bf Lagrangian} if the restriction of $\Omega$ to
$X$ is zero and $\dim(X)=n$, in other words, $X$ is
isotropic of maximal possible dimension.

It is not hard to see that the holomorphically symplectic
form defines the complex structure in a unique way
(\ref{_C_symp_then_HS_Theorem_}). This gives a
way to describe the holomorphically symplectic
structures symplectically, without referring
to the complex structure. We define {\bf a C-symplectic
structure} (\ref{_C_sy_Definition_}) on a real manifold of dimension $4n$
as a closed complex-valued form $\Omega$
such that $\Omega^{n+1}=0$ and $\Omega^n \wedge\bar\Omega^n$
is non-degenerate. Treated in this fashion, the
holomorphic symplectic structure gains the features
we could expect from the symplectic form. In 
\cite{_Soldatenkov_Verbitsky:Moser_}, a C-symplectic version of Moser's
lemma was proven. Just like in case of Moser's lemma,
it was shown that a smooth family $\Omega_t$
of cohomologous C-symplectic structures is trivialized
by an appropriate flow of diffeomorphisms.
However, a significant caveat applies: Moser's
lemma fails unless one requires the vanishing of
the cohomology group $H^1(\calo_{M, I_t})$,
where $I_t$ the family of complex structures on $M$ induced by
the C-symplectic forms $\Omega_t$, and $\calo_{M,I_t}$ is the sheaf of holomorphic functions on $M$ with the complex structure $I_t$.

Weinstein normal neighbourhood theorem is one
of the most immediate and important applications 
of Moser's lemma. This theorem, stated in the
traditional symplectic context, identifies
a tubular neighbourhood of a Lagrangian submanifold
with its neighbourhood in its cotangent space,
equipped with the standard symplectic structure.

This paper is an attempt to apply this reasoning
to C-symplectic manifolds. However, the condition
$H^1(\calo_{M, I_t})=0$ is quite restrictive. It is not always true, and can otherwise be
hard to establish, that $H^1(\calo_{M, I_t})=0$
for a neighbourhood of a holomorphic Lagrangian
submanifold. One exception is when
the holomorphic Lagrangian submanifold
$S\subset M$ is bimeromorphically
contractible, that is, there exists
a holomorphic, bimeromorphic  map 
$M \stackrel \pi \arrow M_0$, bijective on $M \backslash S$, 
which contracts $S$ to a point $s$. 
As follows from Grauert-Riemenschneider theorem,
in this situation the preimage $V:=\pi^{-1}(U)\subset M$ of a Stein neighbourhood $U$ 
of $s\in M_0$ satisfies $H^1(\calo_V)=0$
(\ref{_Gr_Rie+Lagrangian_deforma_Corollary_}).

To establish the existence of a normal form,
we prove that the neighbourhood $V\supset S$ 
can be deformed to a neighbourhood of $S$ 
in its cotangent bundle. For complex manifolds 
this construction, called ``deformation to the
normal cone'', is quite standard. In \ref{_HS_normal_cone_Theorem_},
we show that the deformation to the
normal cone can be performed in the
holomorphically symplectic category.
We construct a smooth family of holomorphic
symplectic manifolds connecting
a neighbourhood of a 
holomorphic Lagrangian submanifold $S\subset M$
to its neighbourhood in $T^* S$.

The content of this paper is split between
Sections \ref{_Moser_Section_} and \ref{_normal_cone_Section_}, 
where the applications of C-symplectic Moser's lemma
are discussed and the deformation to the normal cone 
is constructed, and Section \ref{section-contr}, 
where we show that any bimeromorphically
contractible Lagrangian submanifold is
biholomorphic to the complex projective space.

\subsection{Holomorphic Weinstein theorem}

Recently, Lagrangian subvarieties have attracted
considerable attention in algebraic geometry of
holomorphic symplectic manifolds. One of the reasons is
their relation to the birational geometry: Lagrangian
subvarieties appear as exceptional loci 
of birational contractions of certain type. These are
well-known to be isomorphic to the projective space
whenever smooth (a brief sketch of proof is given in the
beginning of Section \ref{section-contr}). However, the
general, non-algebraic and non-compact case, apparently
has not been studied up to now. In this situation, the
exceptional subvariety of a contraction is defined as
follows.




\hfill

\definition
Let $X\subset Y$ be a complex subvariety in a complex manifold $Y$.
We say that {\bf $X$ can be bimeromorphically contracted}
if there exists a proper morphism of complex varieties
$Y \arrow Y_0$ mapping $X$ to a point, which is an isomorphism 
outside of $X$.

\hfill

Our main result is an analogue of Weinstein tubular
neighbourhood theorem for bimeromorphically contractible
Lagrangian subvarieties, as follows.

\hfill

\theorem\label{_Weinstein_Intro_Theorem_}
Let $(M,I, \Omega)$ be a holomorphically symplectic manifold
(not necessarily compact or quasi-projective), and $E\subset (M,I)$ a compact 
holomorphic Lagrangian submanifold. Assume that $E$ can be bimeromorphically
contracted to a point.  Then $E$ is isomorphic to $\C P^n$.
Moreover,  $E$ has a neighbourhood which is biholomorphically symplectomorphic
to a neighbourhood of the zero section in $T^*\C P^n$.

\hfill

The first statement shall be proved in the next section,
and the proof of the second part shall be given in
Subsection \ref{_deforma_cone_hs_Subsection_}, after
recalling Moser's isotopy lemma and its holomorphic
version in Section 3.

\hfill

\remark Weinstein tubular neighbourhood 
theorem fails when $E$ is a fiber of a holomorphic Lagrangian fibration
on a hyperk\"ahler manifold (say, on an elliptic K3 surface). 
Indeed, the normal bundle $NE$ is trivial, but the elliptic
curve in the elliptic family varies, hence its
neighbourhood cannot be isomorphic to a neighbourhood of
the zero section in the total space of $T^* E= E \times \C$.


\section{Contractible Lagrangian submanifolds}\label{section-contr}


The main result of this section is the following theorem.

\hfill

\theorem\label{_contra_CP^n_Theorem_}
 Let $X$ be a bimeromorphically contractible compact Lagrangian
submanifold of a holomorphic symplectic
manifold $M$ of dimension $2n$, not necessarily
compact or quasiprojective. 
 Then $X$ is isomorphic to $\C P^n$.

 \hfill

This section is devoted to the proof, in fact we shall provide two versions of the argument. These versions have a common beginning which proves that $X$ is Moishezon. Then one may proceed either by adaptation of log-MMP\footnote{This part has been added after the first version of this paper has been written, and is suggested by the subsequent discussions with C. Shramov and Yu. Prokhorov.} following \cite{Kaw}, \cite{Nak}, or by using the positivity properties of the cotangent bundle restricted to a movable curve (see \cite{CP}, \cite{CPT}).

\hfill

To begin with, we recall that in the projective setting,
this result is well-known, see 
e.g. 
\cite{_CMS-B:AdvStud_}, chapter 8.   
The main steps of the proof in the projective case (that
is, assuming the contraction morphism $f: M\to M_0$
projective, as well as $M$ and $M_0$) are 
briefly sketched as follows.

\begin{itemize}

	\item {\bf Step 1:} The contraction morphism is a
          log-extremal contraction with respect to a
          certain boundary divisor $\Delta$
on $M$. More precisely, one takes $H$ ample om $M_0$, then $f^*H$ is big, and writes $f^*H=A+E$
  where $A$ is ample on $M$ and $E$
		effective on $M$. Then $(M, \epsilon E)$ is a klt-pair for a small rational number $\epsilon$, and the contraction is log-extremal for this pair, since $E$ is $f$-negative. By
the log Minimal Model theory, the contraction locus $X$ is
uniruled, i. e. covered by rational curves (see
e.g. \cite{Kaw}).

	\item {\bf Step 2:} According to Z. Ran \cite{Ra},
          Corollary 5.1, any rational curve in a
		holomorphic symplectic manifold $M$ (not necessarily algebraic or compact) deforms in a
          family of dimension at least $2n-2$.

	\item {\bf Step 3:} 
If $C\subset X$ is a rational curve, it is contracted by
$f$, so that any deformation of $C$ inside $M$ must be
contracted as well and therefore lies in $X$. Hence any
rational curve in $X$ deforms in a family of dimension at
least $2n-2$. In particular this is true for a covering
family of minimal rational curves\footnote{I.e. such that
  all of its members passing through a general point are
  irreducible. Informally, $X$ should not covered by
  rational curves splitting off some curves in this
  family.}. Now the results by S. Kebekus 
\cite{_Kebekus_}
imply that this is possible 
only when $X\cong \P^n$.\footnote{See \cite{_CMS-B:AdvStud_}, Theorem 4.2. Since
\cite{_CMS-B:AdvStud_} is known to be not always accurate, we prefer
to rely on  \cite{_Kebekus_}, who proves in Proposition 3.1 that the
base
$H_x$ of a family of minimal rational curves passing through the
general point $x$ and covering the variety $X$ is $\P^{n-1}$, and in
the following subsection 3.2 computes that the total space $U_x$ of
this family
is $\P({\cal O}_{\P^{n-1}}\oplus {\cal O}_{\P^{n-1}}(-1))$, so that
the evaluation map contracts the exceptional section onto $\P^n$.}

\end{itemize}

In this section we give a proof of
\ref{_contra_CP^n_Theorem_}
which does not require any
projectivity assumption. We start with a contraction
criterion by Ancona and Vo Van Tan \cite{AnVo}, and obtain
in particular that $X$ is Moishezon (i.e. bimeromorphic to
projective). Then the results 
by Campana--Peternell--Toma and Campana--Paun \cite{CP}
(formulated in the projective context in \cite{CP}, but
remarked by the authors to be of
bimeromorphic nature) imply uniruledness. Applying Step 2,
we obtain a $2n-2$-dimensional family of minimal rational
curves on $X$. It remains to remark that Kebekus' proof is
also valid when $X$ is Moishezon.

Alternatively (we thank Yuri Prokhorov and Costya Shramov for many indications on this argument), 
once we know that $X$ is Moishezon we can remark that either $X$ has rational curves,
or $X$ is projective (or both) (\cite{Sh}). Then in the first case we repeat Step 2 and Step 3,
and in the second case we make a small adjustment to reduce to the projective 
case.

\hfill

\remark\label{grauert-nbhd} 
As soon as we know that
$X\cong \P^n$, we can apply Grauert's cohomological
criterion (\cite{_Grauert:Modifikationen_}, \S 4, Satz 7,
Corollar) to show that a neighbourhood of $X$ in $M$ is
biholomorphic to a neighbourhood of the zero section in
the cotangent bundle of $\P^n$: see
\cite{_CMS-B:AdvStud_}, Example 8.1 and Lemma 8.2. 

\hfill

Later on in this paper we also prove 
the Weinstein normal form theorem for 
the contractible Lagrangian submanifold
(Subsection \ref{_deforma_cone_hs_Subsection_}), i.e.,
that there exists a biholomorphism taking the symplectic
form on $M$ into the standard symplectic form on the
cotangent bundle.

\subsection{Vo Van Tan's characterization of contractible subvarieties}

Here we recall the main results about the 
bimeromorphically contractible
subvarieties, following \cite{_Grauert:Modifikationen_}
\cite{_VoVanTan:Grauert_conj_}, \cite{AnVo}. All varieties
and analytic spaces are assumed to be connected.

\hfill

\definition
A complex analytic space $X$ is called {\bf  1-convex} 
if it admits a proper holomorphic, bimeromorphic map $p:X \to Y$ to a
Stein space $Y$, and $p_*{\cal O}_X={\cal O}_Y$.

\hfill

\remark
By Remmert Reduction theorem (\cite{_Remmert:Reduction_}), 
a complex variety is holomorphically
convex if and only if it admits a proper, holomorphic map
with connected fibers to a Stein variety. 
The 1-convexity is a stronger property which requires
this map to be bimeromorphic. 

\hfill

\definition
Given a coherent sheaf $A$ on $X$, denote by $L(A)$
the complex analytic space obtained as the relative
spectrum of the sheaf of rings 
$\bigoplus \Sym^k(A)$ over $X$.

\hfill

\remark
When $A$ is a vector bundle, $L(A)$ is the total space 
of the dual bundle $A^*$.

\hfill

\definition
A coherent sheaf $A$ on a compact 
complex analytic space $X$
is called {\bf ample} if for any 
coherent sheaf ${\cal F}$ on $X$ there exists
$k>0$ such that $\Sym^k(A) \otimes {\cal F}$
is globally generated. It is called
{\bf cohomologically positive}
if there exists
$k>0$ such that $H^i(X, \Sym^k(A) \otimes {\cal F})=0$
for all $i>0$, and {\bf weakly positive}
if the zero section of $L(A)$ admits a 1-convex
neighbourhood. 

\hfill

\claim
Suppose that $A$ is a torsion-free coherent sheaf on $X$.
Then $A$ is weakly positive if and only if the zero 
section $X\subset L(A)$ is bimeromorphically contractible.

\hfill



{\bf Proof. Step 1:}
The group
$\C^*$ acts on $L(A)$ by taking any 
compact subvariety which is not contained
in the zero section to a non-compact
family of compact subvarieties.
Since 1-convexity implies that 
all positive-dimensional compact subvarieties of $L(A)$
are contained in a compact subset,
this is impossible, and 
any bimeromorphically contractible
subset of $L(A)$ is contained in the zero section.

\hfill

{\bf Step 2:} The 1-convexity means
that $L(A)$ is equipped with a
bimeromorphic contraction $f$ to a Stein 
variety. Since the zero section $X\subset L(A)$ 
is compact, and a Stein variety has
no compact complex subvarieties, 
the image $f(X)$ is zero-dimensional,
hence it is contracted. As follows
from Step 2, nothing else is contracted.
\endproof

\hfill

Weak positivity has been introduced by Grauert in
\cite{_Grauert:Modifikationen_}, who proved that for line
bundles on compact complex manifolds weak positivity is
equivalent to positivity, i.e. ampleness, and that a
compact complex space carrying a weakly positive vector
bundle is projective algebraic
(\cite{_Grauert:Modifikationen_}, \S 3, Satz 1 and 2).
These results were subsequently extended by Vo Van Tan,
\cite{_VoVanTan:Grauert_conj_}, as follows.

\hfill

\theorem (\cite[Theorem 1]{_VoVanTan:Grauert_conj_})
Let $A$ be a coherent sheaf on a compact irreducible complex space $X$.
Then the three conditions defined above
(ampleness,  weak positivity, cohomological positivity)
are equivalent. \endproof


\hfill

\theorem 
(\cite[Corollary 6]{_VoVanTan:Grauert_conj_}) A
compact irreducible complex space carries a torsion-free
weakly positive coherent sheaf if and only if it is Moishezon.
\endproof

\hfill

\remark\label{contrac-vs-positive} Grauert proved in
\cite[\S 3, Satz 8]{_Grauert:Modifikationen_},  that a
compact submanifold $X$ in an analytic space $M$ with
weakly positive conormal bundle is bimeromorphically
contractible, but observed that the converse is not true,
that is, the conormal bundle of a contractible submanifold
does not have to be weakly positive. He asked whether the
contractibility of $X$ implies that the conormal sheaf of
{\it some} analytic space structure on $X$ is weakly
positive. This question has been answered affirmatively by
Vo Van Tan and Ancona.









More precisely, there is the following contractibility criterion:

\hfill

\theorem (\cite[Corollary 7]{_VoVanTan:Grauert_conj_}, \cite[Corollary 3]{AnVo})\label{contractibility}
Let $M$ be a complex analytic space and  $X\subset M$ a
compact subvariety. Then $X$ is bimeromorphically
contractible if and only if there exists
an ideal sheaf ${\cal J}\subset \calo_M$
such that $\Supp(\calo_M/{\cal J})=X$  and the conormal
sheaf, which is defined as ${\cal J}/{\cal J}^2$,
is ample on $X$. 

\hfill

\remark To say that $X$ carries an ample coherent
torsion-free sheaf is equivalent to saying that there is
an ample coherent sheaf with sheaf-theoretic support equal to $X$.

\hfill

\remark\label{Moishezon-ratcurves} 
In particular, it
follows that a bimeromorphically contractible subvariety
is Moishezon.

\subsection{Contractible subvarieties are isotropic}

Recall the following classical theorem of complex geometry
by Grauert and Riemenschneider \cite{GR}, Satz 2.4 (see
also \cite{T}, Corollary I). 

\hfill

\theorem \label{_Grauert_Riemenschneider_Theorem_}
Let $f:\; U \arrow V$ be a bimeromorphic contraction of
complex spaces with $U$ smooth.
Then  $R^i f_* (K_U)=0$ for $i>0$, where $K_U$ is the
canonical bundle of $U$.

\hfill

\corollary\label{_Gr_Rie+Lagrangian_deforma_Corollary_}
Let $X\subset M$ be a bimeromorphically
contractible subvariety
of a manifold with trivial canonical bundle. Then 
any open neighbourhood of $X$ in $M$ contains
a tubular neighbourhood $U\supset X$ such that 
$H^i(\calo_U)=0$ for any $i>0$.

\hfill

\proof
Let $f:\; M \arrow M_0$ be the bimeromorphic
contraction, mapping $X$ to a point $x\in M_0$.
Consider a Stein neighbourhood $V \ni x$, and
let $U:= f^{-1}(V)$.
By assumption, $K_M = \calo_M$.
Grauert-Riemenschneider theorem
(\ref{_Grauert_Riemenschneider_Theorem_}) implies that $R^i f_*(\calo_U)=0$.
The Grothendieck spectral sequence
with $E_2$-table $H^j(R^i f_*(\calo_U))$
converges to $H^{i+j}(\calo_U)$, giving
$H^k(\calo_U)= H^k(f_* \calo_U) = H^k(\calo_V)=0$
because $V$ is Stein. 
\endproof

\hfill

The fact that a contractible subvariety is Moishezon,
together with Grauert-Riemenschneider theorem, yields the
following result which is well-known in the projective
case (see e. g. \cite{_Kaledin:Poisson_singularities_}): 

\hfill

\proposition Let $M$ be a holomorphic symplectic manifold
and $X$ a contractible subvariety. Then $X$ is isotropic,
that is, the restriction of the holomorphic symplectic
form $\Omega$ to $X$ is zero. 

\hfill

\proof We may assume that $M$ is a neighbourhood of $X$ as
above, so that $H^2(\calo_M)=0$. Let $h: X'\to X\subset M$
be a resolution of singularities of $X$ with $X'$
projective (it exists because $X$ is Moishezon, 
\ref{Moishezon-ratcurves}). The
Dolbeault cohomology class of the $(2,0)$-form
$h^*\overline{\Omega}$ is zero, because the 
Dolbeault cohomology class of 
$\overline{\Omega}$ itself is zero on $M$. Since $X'$ is
smooth projective, by Hodge theory $h^*{\Omega}=0$, so
that $\Omega$ must be zero on $X$. \endproof

\hfill

\subsection{An MMP-style proof of \ref{_contra_CP^n_Theorem_}} 

Let $f: M\to M_0$ be our bimeromorphic contraction. Without loss of generality, we replace $M_0$ by a small Stein neighbourhood of $x=f(X)$. In particular, we may assume that $M$ retracts on $X$. We have established that $X$ is Moishezon.
It is known, see \cite{P} in
dimension 3, \cite{Sh}, \cite{VP} in arbitrary dimension,
that all Moishezon varieties either contain rational curves,
or are projective, or both.

\medskip

{\bf Case 1:} $X$ contains a rational curve. Since $X$ is isotropic, we know that $\dim(X)\leq n$.
By \cite{Ra}, any rational curve in $M$ deforms in a family of dimension at least $2n-2$. Since 
$M_0$ is Stein, all rational curves in $M$ are contained in the contraction locus $X$. We claim they cover $X$ and the dimension of $X$ is $n$. Indeed without loss of generality we may consider a family of ``minimal'', or ``generically unsplit'' curves: this means that all curves through a general point $z$ of the subvariety $Z\subset X$ which they cover, are irreducible.
Then by bend-and-break there is only a finite number of minimal curves through $z$ and any other point of $Z$. Dimension count gives $\dim(Z)\geq n$ and $Z=X$. 
One proceeds to the proof of \ref{_CP^n_main_proof_Theorem_} to finish the argument.

\medskip

{\bf Case 2:} $X$ is projective.

\hfill

\claim Let $M$ be a complex manifold with trivial
canonical bundle admitting a 
holomorphic birational contraction $f:\; M \arrow M_0$,
contracting a submanifold $X\subset M$ to a point $x$. If $X$ is projective, 
then some ample line bundle $L$ on $X$ extends, together with its sections,
to a line bundle $L'$ on $f^{-1}(U)$, where $U$ is a suitable neighbourhood of $x$.

\hfill

\proof 
We choose $U$ in such a way that it is Stein and
$X$ is a deformation retract of $f^{-1}(U)$.
To simplify the notation, we may assume 
that $U=M_0$, and $f^{-1}(U)=M$.
Consider the exponential exact sequence
$0 \arrow \Z_M\stackrel{2\pi\1}\arrow \calo_M \arrow
\calo_M^*\arrow 0$. From this exact sequence (and a
similar sequence on $X$) we obtain the following
diagram, with all rows exact.
\[
\begin{CD}
H^1(M, \calo_M) @>>> H^1(M, \calo_M^*) @>>> H^2(M, \Z) @>>> H^2(M, \calo_M) \\
@VVV @VVV @VVV@VVV \\
H^1(X, \calo_X) @>>> H^1(X, \calo_X^*) @>>> H^2(X, \Z) @>>> H^2(X, \calo_X) 
\end{CD}
\]
The vertical arrow $H^2(M, \Z) \arrow H^2(X, \Z)$ 
in this diagram is an isomorphism,
and $H^2(M, \calo_M)=0$ by Grauert-Riemenschneider 
(\ref{_Grauert_Riemenschneider_Theorem_}).
This implies that any $\eta\in H^2(X, \Z)$
which belongs to the image of $H^1(X, \calo_X^*)$
can be lifted to an element of $H^2(M, \Z)$
which belongs to the image of $H^1(M, \calo_M^*)$.
Therefore, the restriction map is surjective on the set of connected
components of $H^1(X, \calo_X^*)$. Since ampleness is a numerical condition, a connected component of $\Pic(X)$
containing an ample line bundle consists of ample line bundles.  
 
Therefore some ample $L$ on $X$ extends to a line bundle $L'$ on $M$. The proof of \cite{Nak}, Proposition 1.4, shows that for a sufficiently large $k$,  
the restriction map $(L')^{\otimes k}\to L^{\otimes k}$ is surjective on the global sections.
\endproof

\hfill

Hence the contraction $f$ is projective. Now we obtain the uniruledness by the proof of \cite{Kaw}, 
Theorem 2. Indeed, this theorem is a consequence of \cite{Kaw}, Lemma 2, where the crucial step is the 
vanishing 
of higher direct images from \cite{Nak}, Theorem 3.6, valid in the complex-analytic context. Once the uniruledness is established, we again proceed to \ref{_CP^n_main_proof_Theorem_}.

\hfill

\remark\label{analytic-mmp} As noted in the proof of \cite{Nak}, Theorem 3.6, we are in fact dealing with a version of 
log-MMP here.
Namely, to apply log-MMP one needs a boundary divisor negative on $X$. In the projective setting, it is given by decomposing the big divisor $f^*H$, $H$ ample, in the sum of ample and effective. In the complex-analytic setting
one remarks that $(L')^{-1}$ is effective, indeed $f_*((L')^{-1})$ has sections because $M_0$ is Stein. This gives the boundary divisor negative on $X$.

\hfill

In the rest of this section, we give another proof of \ref{_contra_CP^n_Theorem_}, where the uniruledness of $X$ is obtained directly from Ancona-Vo Van Tan criterion with the help of \cite{CP}.

\subsection{Lagrangian contractible submanifolds are uniruled}

Here we apply Vo Van Tan's results to a 
bimeromorphically contractible Lagrangian
submanifold $X$ and prove that it is uniruled
(that is, covered by rational curves). 

\hfill

\lemma\label{_tensor_ample_Lemma_}
 i) A quotient of an ample sheaf is ample;

ii) A tensor power of an 
ample sheaf is again ample;

iii) If $\Supp({\cal A})=X$, the restriction of an ample
sheaf ${\cal A}$ to a curve on $X$, which is not contained
in the singular set of ${\cal A}$, is of strictly positive
degree.

\hfill

\proof The first assertion is true because the $k$-th
symmetric power of a sheaf surjects onto the $k$-th
symmetric power of its quotient, the second is proved
using the representation theory, exactly as in the vector
bundle case treated in \cite{H}, and the third results
from the fact that a globally generated sheaf of positive
rank on a curve has strictly positive degree, unless 
it is trivial.  \endproof

\hfill

\lemma\label{_ample_in_ideal_Lemma_}
Let $X\subset M$ be a smooth submanifold,  and
${\cal J}\subset \calo_M$ an ideal sheaf defining an
analytic subspace of $M$ with support $X$ (that is, 
$\Supp (\calo_M/{\cal J})=X$).  Denote by 
$J_X$ the ideal of $X$. Then there exists
a non-zero map from ${\cal J}/{\cal J}^2$
to $\Sym^r(J_X/J_X^2)= J_X^r/J^{r+1}_X$ for some $r$.


\hfill

\proof
We have ${\cal J}\subset J_X$.
Choose $r$ such that ${\cal J} \subset J_X^r$ and
${\cal J}\subsetneq J_X^{r+1}$.
Then the natural map
${\cal J}\arrow J_X^r/J^{r+1}_X$
is non-zero. However, since
${\cal J}\subset J_X$,
this map vanishes on ${\cal J}^{2}$,
defining a non-zero sheaf morphism
${\cal J}/{\cal J}^{2}\arrow J_X^r/J^{r+1}_X$.
\endproof

\hfill

\corollary If $X$ is a contractible Lagrangian submanifold in a holomorphic symplectic $M$, then for some $r$, $\Sym ^r TM$ contains an ample 
subsheaf.

\hfill

\proof For such an $X$, the holomorphic symplectic form provides an isomorphism between the conormal bundle and the tangent bundle.
Take ${\cal J}$ as in \ref{contractibility} and apply the above lemma.
\endproof

\hfill

In the sequel, we will need the following theorem,
implicit in the work of Campana and Paun \cite{CP} and the preceding work of Campana, Peternell and Toma \cite{CPT}. Define a  {\bf movable curve} $C$ on $X$ as a member of a dominating family of irreducible curves.

\hfill

\theorem\label{_ample_then_uniruled_Theorem_}
Let $X$ be a Moishezon manifold.
Assume that for some $r>0$, 
$\Sym ^r TX$ has a subsheaf ${\cal E}$ such that $deg({\cal E}|_C)>0$ for a movable curve $C$.
Then $X$ is uniruled.

\hfill

\proof Let $\alpha$ be the cohomology class of $C$, then, in the terminology of \cite{CP}, $\mu_{\alpha, max}({\Sym^r TX})>0$.
Then also $\mu_{\alpha, max}(TX)>0$ (\cite{CP}, Theorem 2.9), and the maximal destabilizing sheaf ${\cal F}$ of $TX$ has positive $\alpha$-slope (together with all its quotients). 
If $X$ was projective, we would conclude immediately that $X$ is uniruled, by
\cite{CP}, Theorem 1.1. Otherwise, we let $\pi: X'\to X$ be a projective modification. The class $\pi^*\alpha$ is a movable class on 
$X'$. Indeed, by \cite{BDPP} a movable class $\beta$ on a projective manifold is characterized by $D\beta\geq 0$ for all effective divisor classes $D$.
Setting $\beta=\pi^*\alpha$, we get that this condition is satisfied by the projection formula: indeed $\alpha \pi_*D=0$ when $D$ is contracted by $\pi$, and is non-negative otherwise.
Now in the same way as in Lemma 2.12 of \cite{CP}, the subsheaf  ${\cal F}'\subset TX'$ defined by ${\cal F}'=\pi^*{\cal F}\cap TX'$ is of positive slope with respect to the movable class $\pi^*\alpha$, and we conclude by the main results of \cite{CP} that $X'$, and hence also $X$, is uniruled.
\endproof

\hfill

\corollary A contractible Lagrangian submanifold of a holomorphic symplectic manifold is uniruled. \endproof



\subsection{Uniruled Lagrangian submanifolds}

Following \cite{_CMS-B:AdvStud_,_Hu_Yau:Birational_,_Kebekus_}, 
we use the uniruledness (\ref{_ample_then_uniruled_Theorem_})
to prove the following theorem.

\hfill

\theorem\label{_CP^n_main_proof_Theorem_}
Let $X\subset M$ be a bimeromorphically
contractible Lagrangian submanifold of a holomorphically
symplectic $2n$-dimensional manifold $M$. Then $X$ is biholomorphic to $\C P^n$.

\hfill

\proof By \cite{Ra}, Corollary 5.1, any rational curve $C$ on $M$ deforms in a family of dimension at least $2n-2$. When $C\subset X$, any deformation of $C$ also lies in $X$, because $X$ is contractible. We conclude that any family of rational curves in $X$, 
$\dim(X)=n$, has dimension at least $2n-2$. Such families exist because $X$ is uniruled.

We apply this to a complete family of minimal rational curves covering $X$. ``Minimal'' in this context means that all curves from the family passing through a sufficiently general point $x$ of $X$ are irreducible. It is easy to see by noetherian induction that such minimal rational curves exist. We obtain a family of irreducible curves through $x$ of dimension at least $n-1$.
These curves must cover $X$, by bend-and-break (yielding that there are only finitely many curves connecting $x$ to another  
point $y$) and by compacity of families of cycles on Moishezon manifolds. Now we are in the situation studied by
Kebekus for projective manifolds, see \cite{_Kebekus_}. The arguments carry over verbatim to the Moishezon case and prove that $X$ is isomorphic to the projective space.
\endproof


\section{Moser's lemma and normal form of Lagrangian submanifolds}
\label{_Moser_Section_}


\subsection{Moser's lemma  for holomorphic symplectic and C-symp\-lectic
structures}

The complex structure on a holomorphically symplectic manifold
is uniquely determined by the holomorphic symplectic
form. This was the starting point of 
 \cite{_Soldatenkov_Verbitsky:Moser_}.
In \cite{_Soldatenkov_Verbitsky:Moser_} the authors obtain a version of 
holomorphically symplectic Moser's lemma 
which accounts for the deformations of complex structure
as well as for the deformations of the holomorphic
symplectic form.

\hfill

Moser's lemma is the following fundamental result of symplectic 
geometry.

\hfill

\theorem
Let $\omega_t$, $t\in [0,1]$ be a smooth family of symplectic
structures on a compact manifold $M$.
Assume that the cohomology class $[\omega_t]\in H^2(M)$ is constant in $t$.
Then there exists a smooth family $\Psi_t\in \Diff_0(M)$
of diffeomorphisms such that $\Psi_t^*\omega_0=\omega_t$.

\proof \cite{_McD-Sal-intro_}.
 \endproof

 \hfill

\remark The proof proceeds by constructing a vector field and integrating it to a flow of diffeomorphisms.
It is for integration that one needs compacity. In particular, one can obtain a similar statement replacing a compact manifold $M$ by a neighbourhood of a compact submanifold $X\subset M$.

\hfill

For a holomorphically symplectic version of this lemma,
we recall the notion of C-symplectic structures. 
A C-symplectic structure is a holomorphically symplectic
form understood abstractly, that is, without fixing
a complex structure (which is uniquely determined by 
the holomorphically symplectic form nevertheless).

\hfill

\definition\label{_C_sy_Definition_}
 (\cite{_BDV:sections_})
Let $M$ be a smooth $4n$-dimensional manifold.
A closed complex-valued form $\Omega$ on $M$
is called {\bf C-symplectic} if $\Omega^{n+1}=0$
and $\Omega^{n}\wedge \bar\Omega^n$ is a non-degenerate
volume form.

\hfill

\theorem\label{_C_symp_then_HS_Theorem_}
Let $\Omega\in \Lambda^2(M,\C)$ be a C-symplectic form,
and $T^{0,1}_\Omega(M):=\ker \Omega$, where
\[
 \ker \Omega:= \{ v\in TM\otimes \C \ \ |\ \  \Omega\cntrct v=0\}.
\]
Then $T^{0,1}_\Omega(M)\oplus 
\overline{T^{0,1}_\Omega(M)} = TM\otimes_\R \C$, hence
the sub-bundle $T^{0,1}_\Omega(M)$ defines an 
almost complex structure $I_\Omega$ on $M$. Moreover, $I_\Omega$ is integrable,  and
$\Omega$ is holomorphically symplectic on $(M, I_\Omega)$.

\proof \cite[Proposition 2.12]{_BDV:sections_},
\cite[Theorem 3.5]{_Verbitsky:Degenerate_}.
\endproof

\hfill

Now we can state the C-symplectic version of Moser's lemma.

\hfill

\theorem 
Let $(M, I_t, \Omega_t)$, $t\in [0,1]$ 
be a family of C-symplectic
forms on a compact manifold, with corresponding complex
structures $I_t$. Assume that the cohomology class
$[\Omega_t]\in H^2(M,\C)$ is constant, and 
$H^{0,1}(M,I_t)=0$, where $H^{0,1}(M,I_t)= H^1(M, \calo_{(M, I_t)})$
is the first cohomology of the sheaf of holomorphic functions.
Then there exists a smooth family of diffeomorphisms 
$V_t\in \Diff_0(M)$, such that $V_t^*\Omega_0=\Omega_t$.

\proof \cite[Theorem 2.5 (2)]{_Soldatenkov_Verbitsky:Moser_}. \endproof

\hfill

For non-compact manifolds, a slightly more cumbersome 
version of Moser's lemma can be stated.

\hfill

\theorem \label{_Moser_noncompact_Theorem_}
Let $\pi\colon {\cal X} \to \Delta$ be
a smooth family of holomorphic
symplectic manifolds (not necessarily compact)
over the unit disc, locally trivial as a family of $C^\infty$ manifolds.
Denote by ${\cal X}_t= \pi^{-1}(t)$ its fiber, and
let $\Omega_t\in H^0({\cal X}_t, \Omega^2_{{\cal X}_t})$ be its holomorphic
symplectic form, smoothly depending on $t$. 
Using the $C^\infty$ trivialization to identify
cohomology groups of the fibres, assume that the cohomology
class of $\Omega_t$ does not depend on $t\in \Delta$,
and $H^1({\cal X}_t, \calo_{{\cal X}_t})=0$. 
Let $K \subset {\cal X}_{t_0}$ be a compact subset. 
Then there exists an open neighbourhood
$U\subset \Delta$  of $t_0\in \Delta$, and an open subset 
${\tilde U} \subset \pi^{-1}(U)$, with $K\subset {\tilde U}$,
with the following property. The set
${\tilde U}$ is locally trivially
fibered over $U$, with all fibres ${\tilde U} \cap \pi^{-1}(t)$,
$t\in U$ isomorphic as holomorphically symplectic
manifolds.

\proof \cite[Theorem 2.5 (1)]{_Soldatenkov_Verbitsky:Moser_}. \endproof

\hfill

We will use an equivalent version of this result.

\hfill

\theorem \label{_Moser_noncompact_family_compacts_Theorem_}
Let $\pi\colon {\cal X} \to \Delta$ be
a smooth family of holomorphic
symplectic manifolds (not necessarily compact)
over the unit disc, trivial as a family of $C^\infty$ manifolds.
Denote by ${\cal X}_t= \pi^{-1}(t)$ its fiber, and
let $\Omega_t\in H^0({\cal X}_t, \Omega^2_{{\cal X}_t})$ be its holomorphic
symplectic form, smoothly depending on $t$. 
Using the $C^\infty$ trivialization to identify
cohomology groups of the fibres, assume that the cohomology
class of $\Omega_t$ does not depend on $t\in \Delta$,
and $H^1({\cal X}_t, \calo_{{\cal X}_t})=0$. 
Let $K_t \subset {\cal X}_{t}$ be a locally trivial
family of compact subsets, and ${\cal K}\subset{\cal X}$ their union. 
Then there exists an open neighbourhood
$U\subset \Delta$  of $t_0\in \Delta$, and an open subset 
${\tilde U} \subset \pi^{-1}(U)$, with $K_t\subset {\tilde U}$,
with the following property. The set
${\tilde U}$ is locally trivially
fibered over $U$, with all fibres ${\tilde U} \cap \pi^{-1}(t)$,
$t\in U$ isomorphic as holomorphically symplectic
manifolds.

\hfill

\proof
Chose $U\subset \Delta$ and $\tilde U\subset {\cal X}$
as in \ref{_Moser_noncompact_Theorem_}. Denote by $\Delta_\lambda$
the disc of radius $\lambda$ centered in 0. 
We need to show that 
${\cal K} \cap\pi^{-1}(\Delta_\lambda) \subset \tilde U$
for $\lambda$ sifficiently small.
Introduce a Riemannian metric on ${\cal X}$. Since $\tilde U\subset {\cal X}$
is open, it contains an $\epsilon$-neighbourhood of $K_0$, for a sufficiently
small $\epsilon >0$.
Choose $\lambda$ in such a way that $K_t$ belongs
to an $\epsilon$-neighbourhood $K_0(\epsilon)$
of $K_0$ for all $t\in \Delta_\lambda$;
this is possible to do, because the Hausdorff distance  $d_H(K_t, K_0)$
converges to 0 as $t$ converges to 0. For such $\lambda$, clearly,
one has 
\[ 
{\cal K} \cap\pi^{-1}(\Delta_\lambda)\subset K_0(\epsilon)\subset \tilde U.
\]
\endproof

\subsection{Moser's lemma  for a neighbourhood of a 
holomorphic Lagrangian submanifold}

We will apply \ref{_Moser_noncompact_family_compacts_Theorem_} 
when $K$ is bimeromorphically
contractible Lagrangian submanifolds. We need the following
version of Moser's lemma.

\hfill

\lemma\label{_lagra_Moser_Lemma_}
Let $(M,I_t, \Omega_t)$, $t\in [0,1]$ be 
a smooth family of C-symplectic manifolds
(not necessarily compact),
with all $\Omega_t$ exact, and $E_t\subset (M,I_t)$ be a family of compact
holomorphic Lagrangian subvarieties. Assume that 
$H^{0,1}(M, I_t)=0$.
Then $E_t$ have open neighbourhoods $U_t$
in $M$ such that  $(U_t, I_t, \Omega_t, E_t)$
is trivialized by a flow of diffeomorphisms.

\hfill

\proof
Follows immediately from
\ref{_Moser_noncompact_family_compacts_Theorem_}.
\endproof

\hfill
 


\hfill










In the situation described by 
\ref{_Gr_Rie+Lagrangian_deforma_Corollary_},
any deformation of a neighbourhood of
a holomorphic Lagrangian subvariety is locally
trivial by \ref{_lagra_Moser_Lemma_}, because the relevant
cohomology groups vanish. This brings

\hfill

\corollary\label{_trivial_defo_CP^n_Corollary_}
Let $(M,I, \Omega)$ be a holomorphically symplectic manifold
(not necessarily compact), and $E\subset (M,I)$ a compact 
holomorphic Lagrangian subvariety biholomorphic to $\C P^n$. 
Assume that a neighbourhood of $E$ can be smoothly deformed
to a neighbourhood of the zero section in $T^* E$
as a C-symplectic manifold. 
Then $E$ has a neighbourhood which is isomorphic
to a neighbourhood of $E$ in $T^*E$ as a holomorphically
symplectic manifold.

\hfill

\proof 
We apply \ref{_lagra_Moser_Lemma_} to the deformation $(U_t, I_t, \Omega_t, E_t)$ connecting a neighbourhood $U_0$ of $E$ 
in $M$ to a neighbourhood $U_1$ of $E$ in $T^*E$. Note that since $E\cong \P^n$, all $E_t$ are isomorphic to $\P^n$ 
as well. Therefore they are all bimeromorphically contractible inside $U_t$: indeed
$N^* E_t\cong TE_t$ is isomorphic to the cotangent bundle of $\P^n$ and hence ample, and we conclude by \ref{contrac-vs-positive}. 
Hence $H^i(\calo_{U_t}), i>0,$ of a tubular
neighbourhood of $E_t$ in $(M, I_t)$ vanish by \ref{_Gr_Rie+Lagrangian_deforma_Corollary_}, and \ref{_lagra_Moser_Lemma_} applies.
\endproof

\hfill

\remark The exactness condition on the symplectic forms is satisfied automatically in a sufficiently small neighbourhood $U_t$ of $E_t$, because $U_t$ retracts on $E_t$.


\section{Deformation to the normal cone of a Lagrangian submanifold}
\label{_normal_cone_Section_}


To apply \ref{_trivial_defo_CP^n_Corollary_}, we need 
to deform a neighbourhood of a contractible Lagrangian submanifold
to its neighbourhood in its cotangent space. This construction is well known
in algebraic geometry under the name ``deformation
to the normal cone''; however, we have to make sure
it works in the holomorphic symplectic category.

\subsection{Deformation to the normal cone: an introduction}

Let $X\subset M$ be a complex submanifold in a complex manifold.
{\bf Deformation to the normal cone} is a holomorphic
deformation of a neighbourhood of $X\subset M$ over the disk
such that its central fiber is the total space of the
normal bundle $NX$, and the rest of the fibers are $M$.
It is obtained as follows.

Let $X \subset M$ be a complex subvariety.
Consider a product $M_1:=M \times \Delta$ of $M$ with the
disk $\Delta$, and let $\tilde M_1$ be the blow-up 
of $M_1$ in $X\times \{0\}$. Denote by
$\tilde \pi_1:\; \tilde {\cal X}\arrow \Delta$ the blow-down
composed with the projection. The preimage $\tilde \pi_1^{-1}(0)$
is a union of two irreducible components, the proper 
preimage of $M\times \{0\}$, denoted $D_1$, 
and the blow-up divisor, denoted $D_2$.

\hfill

\definition
The {\bf deformation to the normal cone}
is the complement $\tilde M :=\tilde {M_1}\backslash D_1$.
considered as a fibration over the disk $\Delta$.

\hfill

Clearly, the central fiber of the natural projection
$\tilde M \arrow \Delta$ is $D_2\backslash (D_1\cap D_2)$.

\hfill

\claim
When $X$ is smooth, and $M$ is a tubular
neighbourhood of $X$ in $M$, the complement 
$D_2\backslash (D_1\cap D_2)$ is naturally isomorphic to $NX$,
and the ``deformation to the normal cone''
family $\tilde \pi:\; \tilde M \arrow \Delta$
is locally trivial in the smooth category.

\hfill

\pstep
The blow-up divisor 
$E= {\Bbb P} N_{M_1} X= {\Bbb P}(NX\oplus \calo_X)$, 
and its intersection with $D_1$
is the set of all $l \in {\Bbb P}(NX\oplus \calo_X)$
tangent to $M \times \{0\}$. We identify this intersection 
with ${\Bbb  P}(NX)$. This gives an isomorphism
$D_2\backslash (D_1\cap D_2) ={\Bbb P}(NX\oplus \calo_X)
\backslash {\Bbb  P}(NX)= \Tot(NX)$.

\hfill

{\bf Step 2:}
Now, the tubular neighbourhood of $X\subset M$ is
diffeomorphic to $\Tot(NX)$,  hence all fibers
of $\tilde \pi:\; \tilde M \arrow \Delta$ are
diffeomorphic.
\endproof

\subsection{Deformation to the normal cone in holomorphic
  symplectic category}
\label{_deforma_cone_hs_Subsection_}

The following theorem, together with 
\ref{_trivial_defo_CP^n_Corollary_}, 
implies the holomorphic version of the Weinstein normal neighbourhood theorem
(\ref{_Weinstein_Intro_Theorem_}).

\hfill

\theorem\label{_HS_normal_cone_Theorem_}
Let $(M,\Omega)$ be a holomorphically symplectic manifold, and
$X\subset M$ a 
holomorphically Lagrangian submanifold. 
Then there exists a smooth, holomorphic deformation of
a neighbourhood of $X$ in $M$ over the disk $\Delta$, such that its
central fiber is biholomorphic to a neighbourhood of $X$ in
$T^*X$, the rest of the fibers are biholomorphic to a
neighbourhood of $X$ in $M$, and the holomorphic symplectic
form on $T^* X$ can be smoothly extended to the holomorphic
symplectic form on the rest of the fibers.

\hfill

\remark
Note that this deformation in complex analytic 
category is already constructed: it is the ``deformation
to the normal cone'' family. However, to apply the C-symplectic
Moser lemma, we need to have a smooth family of holomorphically
symplectic forms on its fibers.

\hfill

{\bf Proof of \ref{_HS_normal_cone_Theorem_}. Step 1:}
Let $\tilde M \stackrel {\tilde \pi}\arrow \Delta$
be the deformation to the normal cone family,
and $t$ the coordinate on $\Delta$.
Locally in $X$ we can write $X$ by a
system of holomorphic equations
$q_1=q_2=...=q_n=0$, and the
holomorphically symplectic
form as $\Omega=\sum_{i=1}^n dp_i \wedge dq_i$.
We are going to prove that $t^{-1}\Omega$
is extended to a non-degenerate form 
on the central fiber of $\tilde \pi$.

Since this extension is unique, if would
suffice to prove that $t^{-1}\Omega$
can be extended to the central fiber 
locally in $X$, and to check that
it is non-degenerate. 

\hfill

{\bf Step 2:} 
The coordinates on the central
fiber of the deformation to the normal cone family
$\tilde M\stackrel {\tilde \pi}\arrow \Delta$
are given by $p_1, ..., p_n, \tilde q_1, ..., \tilde q_n$.
Trivializing the neighbourhood of $x\times \Delta\in X\times \Delta$
along $\Delta$ in the usual way, we write 
$\tilde q_i = t^{-1} q_i$: this is the standard 
way to write coordinates on the blow-up. 

Locally, we can always find a $\theta\in \Omega^1 M$ such that
$d\theta = \Omega$ and $\theta \restrict X=0$.
For example, we can take $\theta = \sum_i q_i dp_i$ 
in the above coordinates. 
Writing $\theta$ in the coordinates $\tilde q_i, p_i$,
we get $\theta = t \sum_i \tilde q_i dp_i$.
Then $d(t^{-1}\theta)= \sum_i d \tilde q_i \wedge dp_i$
is a holomorphic form on $\tilde M$. Restricted
to the fibers $\pi^{-1}(u)$ of the projection
$\tilde M\stackrel {\tilde \pi}\arrow \Delta$,
this form is equal to $t^{-1}\Omega + d(t^{-1}) \wedge \theta$;
the second term restricted to $\pi^{-1}(u)$ vanishes,
hence the restrictions of $d(t^{-1}\theta)$ extend
$t^{-1}\Omega$ holomorphically to $\tilde M$.

\hfill

{\bf Step 3:} It remains to show that 
$d(t^{-1}\theta)$ is non-degenerate on
the central fiber of $\tilde M \stackrel {\tilde \pi}\arrow \Delta$.
Writing $\Omega = \sum_{i=1}^n dp_i \wedge dq_i$
as above and passing to the coordinates $q_i = t \tilde q_i$,
we obtain $\Omega = \sum_{i=1}^n t dp_i \wedge d\tilde q_i
+ \sum_{i=1}^n dp_i \wedge \tilde q_i dt$.
Since the last term vanishes on the fibers,
the form $t^{-1} \Omega = \sum_{i=1}^n dp_i \wedge d\tilde q_i$
is smooth, non-degenerate on the central fiber, and
equal to $t^{-1}\Omega$ on the general fibers.
This proves \ref{_HS_normal_cone_Theorem_}.
\endproof

\hfill

{\bf Acknowledgements:} We are grateful to Fr\'ed\'eric Campana, Andreas H\"oring, 
Dmitry Kaledin, Yuri Prokhorov, Costya Shramov and Andrey Soldatenkov for their insightful discussions,
and to Arnaud Beauville for his valuable email communication.

\hfill

{

}
{\small
\noindent {\sc Ekaterina Amerik\\
{\sc Laboratory of Algebraic Geometry,\\
National Research University HSE,\\
Department of Mathematics, 6 Usacheva Str., Moscow, Russia,}\\
\tt  Ekaterina.Amerik@gmail.com}, also: \\
{\sc Universit\'e Paris-11,\\
Laboratoire de Math\'ematiques,\\
Campus d'Orsay, B\^atiment 425, 91405 Orsay, France}

\hfill

\noindent 
{\sc Misha Verbitsky\\
Instituto Nacional de Matem\'atica Pura e
              Aplicada (IMPA) \\ Estrada Dona Castorina, 110\\
Jardim Bot\^anico, CEP 22460-320\\
Rio de Janeiro, RJ - Brasil\\
{\tt  verbit@impa.br} \\
also:\\
{\sc Laboratory of Algebraic Geometry,\\
National Research University HSE,\\
Department of Mathematics, 6 Usacheva Str.,\\ Moscow, Russia}}.

\end{document}